\documentclass[journal]{IEEEtran}
\newif\ifpdf
\ifx\pdfoutput\undefined
	\pdffalse
\else
	\pdfoutput=1
	\pdftrue
\fi
\ifpdf
	\usepackage{graphicx,graphics,latexsym,verbatim,epsfig}
	\usepackage{graphicx}
	\DeclareGraphicsExtensions{.pdf,.jpg,.png}
\else
	\usepackage{graphicx}
	\DeclareGraphicsExtensions{.eps}
\fi
\let\ifpdf\relax

\makeatletter
\def\maxwidth{
  \ifdim\Gin@nat@width>\linewidth
    \linewidth
  \else
    \Gin@nat@width
  \fi
}
\newdimen\mymathindent
    {\@beginparpenalty\predisplaypenalty
     \@endparpenalty\postdisplaypenalty
     \refstepcounter{equation}%
     \trivlist \item[]\leavevmode
       \hb@xt@\linewidth\bgroup $\m@th
         \displaystyle
         \hskip\mymathindent}%
        {$\hfil 
         \displaywidth\linewidth\hbox{\@eqnnum}%
       \egroup
     \endtrivlist}
\makeatother

\usepackage{placeins}
\usepackage[cmex10]{amsmath}
\usepackage{amsfonts}
\usepackage{amsthm}
\usepackage{amssymb}
\usepackage{mathtools}
\usepackage{multirow}
\usepackage{enumerate}
\usepackage[caption=false,font=footnotesize]{subfig}

\newtheorem{theorem}{Theorem}

\newtheorem{corollary}{Corollary}
\newtheorem{definition}{Definition}
\newtheorem{lemma}{Lemma}

\newtheorem{remark}{Remark}

\newcommand{\refe}[1]{(\ref{#1})}

\newcommand{\Mzero}{\ensuremath \left.0\mkern-6.3mu \raisebox{1.0pt}{\textit{\scriptsize /}}\right.}

\newcommand{\Exp}[1]{\mathrm{e}^{#1}}

\begin{document}
\title{A New Notion of Effective Resistance for Directed Graphs---Part I: Definition and Properties}
\author{George~Forrest~Young,~\IEEEmembership{Student~Member,~IEEE,
} Luca~Scardovi,~\IEEEmembership{Member,~IEEE,} and~Naomi~Ehrich~Leonard,~\IEEEmembership{Fellow,~IEEE}%
\thanks{This research was supported in part by AFOSR grant FA9550-07-1-0-0528, ONR grant N00014-09-1-1074, ARO grant W911NG-11-1-0385 and the Natural Sciences and Engineering Research Council (NSERC) of Canada.}%
\thanks{G. F. Young and N. E. Leonard are with the Department of Mechanical and Aerospace Engineering, Princeton University, Princeton, NJ 08544, USA; e-mail: \texttt{gfyoung@princeton.edu}; \texttt{naomi@princeton.edu}.}%
\thanks{L. Scardovi is with the Department of Electrical and Computer Engineering, University of Toronto, Toronto, ON, M5S 3G4, Canada; e-mail: \texttt{scardovi@scg.utoronto.ca}.}}

\date{\today}
\maketitle

\begin{abstract}
The graphical notion of effective resistance has found wide-ranging applications in many areas of pure mathematics, applied mathematics and control theory. By the nature of its construction, effective resistance can only be computed in undirected graphs and yet in several areas of its application, directed graphs arise as naturally (or more naturally) than undirected ones. In part I of this work, we propose a generalization of effective resistance to directed graphs that preserves its control-theoretic properties in relation to consensus-type dynamics. We proceed to analyze the dependence of our algebraic definition on the structural properties of the graph and the relationship between our construction and a graphical distance. The results make possible the calculation of effective resistance between any two nodes in any directed graph and provide a solid foundation for the application of effective resistance to problems involving directed graphs.
\end{abstract}


\section{Introduction}\label{sec:intro}
\IEEEPARstart{T}{he} concept of \emph{effective resistance} has been used in relation to graphs for some time \cite{Klein93}. This concept stems from considering a graph, consisting of a set of nodes connected by weighted edges, to represent a network of resistors (one resistor corresponding to each edge) with resistances equal to the inverse of the corresponding edge weights. Then, the effective resistance between nodes $k$ and $j$, denoted $r_{j,k}$, can be found by the resistance offered by the network when a voltage source is connected between these two nodes. One of the useful properties of the effective resistance is that it defines a distance function on a graph that takes into account all paths between two nodes, not just the shortest path \cite{Klein93}. This allows the effective resistance to be used in place of the shortest-path distance to analyze problems involving random motion, percolation and flows over networks.

Effective resistance has proven to have a number of interpretations and applications over a wide variety of fields. One of the earliest interpretations was in the study of random walks and Markov chains on networks \cite{Doyle84, Tetali91, Chandra96, Palacios01}, where the effective resistance between a pair of nodes was related to expected commute, cover and hitting times and the probabilities of a random walk reaching a node or traversing an edge. More direct applications have arisen in the study of power dissipation and time delays in electrical networks \cite{Ghosh08}. In addition, the effective resistance has been shown to have combinatorial interpretations, relating to spanning trees and forests \cite{Shapiro87} as well as the number of nodes and edges in the graph \cite{Klein93}. Following the work of Klein and Randi\'{c} \cite{Klein93}, there has been a substantial literature investigating the use of effective resistance and the Kirchhoff index in the study of molecular graphs \cite{Bonchev94, Gutman96, Zhu96, Lukovits99, Klein01}. 

More recently, effective resistance has arisen in control theory, in the study of control, estimation and synchronization over networks. Barooah and Hespanha described in \cite{Barooah06} how effective resistance can be used to measure the performance of collective formation control, rendezvous and estimation problems. They developed the theory of estimation from relative measurements further in \cite{Barooah07, Barooah08}. A number of authors have demonstrated the use of effective resistance in investigating the problem of selecting leaders to maximize agreement or coherence in a network \cite{Patterson10, Clark11, Fardad11, Fitch13}. D\"{o}rfler and Bullo used effective resistance in their study of synchronization of oscillators in power networks \cite{Dorfler10}, and subsequently developed a theoretical analysis involving resistance for a graph reduction technique \cite{Dorfler13}. We have also used the concept of effective resistance to measure the robustness to noise of linear consensus over networks \cite{Young10, Young11} as well as the performance of nodes in networks of stochastic decision-makers \cite{Poulakakis12}.

By the nature of its definition, effective resistance is restricted to undirected graphs, and in many applications, including the study of molecular graphs and electrical networks, it is natural to focus solely on undirected graphs. However, in many other applications, including random walks and networked control systems, directed graphs arise just as naturally as undirected ones (which can be thought of as special cases of directed graphs in which every edge also exists in the opposite direction and with equal weight). For example, if the nodes in a graph represent agents and the edges interactions, directed edges result from interactions in which one agent reacts to another but the second either has no reaction to the first or reacts with a different strength or rule. Only in highly constrained circumstances would such a network result in an undirected graph. 

Accordingly, it would be particularly useful if the concept of effective resistance could be extended to apply to any directed or undirected graph, so that analysis that is currently only applicable to undirected graphs could be applied in the more general case. Indeed in \cite{Barooah06, Barooah07, Barooah08}, the authors investigated directed measurement graphs, but assumed undirected communication in order to analyze their systems using effective resistance. Similarly, in \cite{Young10, Young11, Poulakakis12}, we began our investigations using directed graphs and then specialized to undirected graphs when we used effective resistance.

In this work we propose a generalized definition of effective resistance for any graph, constructed in such a way that it preserves the connection between effective resistance and networked control and decision systems \cite{Young10, Young11, Poulakakis12}. This new definition produces a well-defined pairwise property of nodes that depends only on the connections between the nodes. Although it is no longer always a metric on the nodes of a graph, our notion of effective resistance does allow for the construction of a resistance-based metric for any graph. This is in contrast with a (perhaps) more intuitive generalization based on the use of pseudoinverses, which does not yield a resistance-based metric in the general case.  Further, this suggests that our construction should prove to be useful for applications other than those we have presented here. In the companion paper we explore some of the implications of our new approach by computing effective resistances in several canonical directed graphs.

This paper is organized as follows. In Section \ref{sec:back} we provide an overview of the notation and definitions used throughout the paper. In Section \ref{sec:def} we define our extended notion of effective resistance. In Section \ref{sec:prop} we analyze some basic properties of our definition, including its well-definedness, its dependence on connections between nodes and its relationship to a metric. We conclude and offer some final thoughts in Section \ref{sec:conc}.

\section{Background and notation}\label{sec:back}
In this section we provide a summary of the notation used throughout this two-part paper and then explain some of the basic definitions and terminology of graph theory, as it applies to this work. There are competing definitions for many of the basic concepts of graph theory, mainly due to varying scope amongst authors. For our purposes, we restrict our attention to directed graphs as they may arise in control theory, and so our definitions mostly follow \cite{OlfatiSaber04}, with some taken from \cite{BangJensen10}. Two notable exceptions are our definition of connectivity, which falls between the standard notions of weak and strong connectivity but is more applicable to control over graphs \cite{Young10, Scardovi09}, and our definition of \emph{connections} in directed graphs, which have interesting parallels to paths in undirected graphs.

Throughout this paper we will represent matrices using capital letters. A matrix $M$ can also be written $\left[m_{i,j}\right]$, where $m_{i,j}$ denotes the scalar entry in the $\left(i,j\right)\text{th}$ position. The identity matrix in $\mathbb{R}^{n\times n}$ will be denoted by $I_n$, and a zero matrix (with dimensions inferred from context) will be denoted by $\Mzero$. Entries of a matrix whose values are not relevant may be denoted by $\ast$. We will represent vectors using bold, lower case letters and their scalar entries with the same (non-bold) lower case letter with a single subscript. In particular, we use $\mathbf{e}^{(k)}_n$ to denote the $k\text{th}$ standard basis vector of $\mathbb{R}^n$. That is, $\mathbf{e}^{(k)}_n$ contains a zero in every position except the $k\text{th}$ position, which is a $1$. In addition, we use $\mathbf{1}_n$ to denote the vector in $\mathbb{R}^n$ containing a $1$ in every entry and $\mathbf{0}$ to denote the zero vector (with size inferred from context). We will use $\text{diag}^{(k)}\!\left(\mathbf{v}\right)$ to denote a $k$-diagonal matrix, with the entries of $\mathbf{v}$ along the $k\text{th}$ diagonal and zeros elsewhere (and the dimensions inferred from the length of $\mathbf{v}$ and $k$). A diagonal matrix will be denoted by $\text{diag}\!\left(\mathbf{v}\right)$.

A \emph{graph} (also \emph{directed graph} or \emph{digraph}) $\mathcal{G}$ consists of the triple $\left(\mathcal{V}, \mathcal{E}, A \right)$, where $\mathcal{V} = \left\{1, 2, \ldots, N \right\}$ is the set of nodes, $\mathcal{E} \subseteq \mathcal{V}\times\mathcal{V}$ is the set of edges and $A \in \mathbb{R}^{N\times N}$ is a weighted adjacency matrix with non-negative entries $a_{i,j}$. Each $a_{i,j}$ will be positive if and only if $\left( i,j \right) \in \mathcal{E}$, otherwise $a_{i,j} = 0$. Note that by viewing $\mathcal{E}$ as a subset of $\mathcal{V}\times\mathcal{V}$, $\mathcal{G}$ can contain at most one edge between any ordered pair of nodes. In addition, we restrict our attention to graphs which do not contain any self-cycles (edges connecting a node to itself).

The graph $\mathcal{G}$ is said to be \emph{undirected} if $\left(i,j\right) \in \mathcal{E}$ implies $\left(j,i\right) \in \mathcal{E}$ and $a_{i,j} = a_{j,i}$. Thus, a graph will be undirected if and only if its adjacency matrix is symmetric. We use the term \emph{undirected edge} to refer to a pair of edges between two nodes (one in each direction), with equal weights.

A graph can be drawn by representing each node with a distinct, non-overlapping circle and representing each edge $(i,j)$ by a line joining the circles for nodes $i$ and $j$. The direction of the edge is indicated by adding an arrow to the line pointing to the circle for node $j$. The edge weight can be written adjacent to the line representing an edge. If no weight is written next to an edge, it is assumed that the edge weight is $1$. An undirected edge (corresponding to $(i,j) \in \mathcal{E}$, $(j,i) \in \mathcal{E}$ and $a_{i,j} = a_{j,i}$) can be represented by a single line either without any arrows or with arrows in both directions, and with the single edge weight written next to the line. 

The \emph{out-degree} (respectively \emph{in-degree}) of node $k$ is defined as $d_k^{\text{\emph{out}}} = \sum_{j=1}^N{a_{k,j}}$ (respectively $d_k^{\text{\emph{in}}} = \sum_{j=1}^N{a_{j,k}}$). A graph is said to be \emph{balanced} if for every node, the out-degree and in-degree are equal. For balanced graphs (including all undirected graphs), the term \emph{degree} is used to refer to both the out-degree and in-degree.

$\mathcal{G}$ has an associated \emph{Laplacian} matrix $L$, defined by $L = D - A$, where $D$ is the diagonal matrix of node out-degrees, that is $D = \mathrm{diag}\left(d_1^{\text{\emph{out}}}, d_2^{\text{\emph{out}}}, \ldots, d_N^{\text{\emph{out}}}\right)$. The row sums of the Laplacian matrix are zero, that is $L \mathbf{1}_N = \mathbf{0}$. Thus $0$ is always an eigenvalue of $L$ with corresponding eigenvector $\mathbf{1}_N$. It can be shown that all the eigenvalues of $L$ are either $0$ or have positive real part \cite{Agaev05}. A graph will be undirected if and only if its Laplacian matrix is symmetric, and then all the eigenvalues of $L$ will be real and non-negative.

The \emph{set of neighbors} of node $k$, denoted $\mathcal{N}_k$, is the set of nodes $j$ for which the edge $(k,j) \in \mathcal{E}$.

A \emph{path} in $\mathcal{G}$ is a (finite) sequence of nodes such that each node is a neighbor of the previous one. A path is called \emph{simple} if no internal nodes (i.e. other than the initial or final nodes) are repeated. The \emph{length} of a path is the number of edges traversed. Thus a single node is considered to be a path of length $0$. A directed (respectively, undirected) \emph{path graph} on $N$ nodes is a graph containing exactly $N - 1$ directed (respectively, undirected) edges and which admits a path of length $N-1$ containing every node in the graph.

A \emph{cycle} in $\mathcal{G}$ is a non-trivial closed path. That is, a cycle is a path of length greater than zero in which the initial and final nodes are the same. A \emph{simple cycle} is a non-trivial closed simple path. Since we are only considering graphs which do not contain self-cycles, the minimum length of a cycle is two. A directed (respectively, undirected) \emph{cycle graph} on $N$ nodes is a graph containing exactly $N$ directed (respectively, undirected) edges and which admits a simple cycle of length $N$ containing every node in the graph. Note however, that while a directed cycle graph on $2$ nodes is possible, the minimum number of nodes in an undirected cycle graph is $3$.

\begin{figure}
\centering
\subfloat{{\footnotesize(a)}\hspace{-0.2cm}\includegraphics[width=3.8cm]{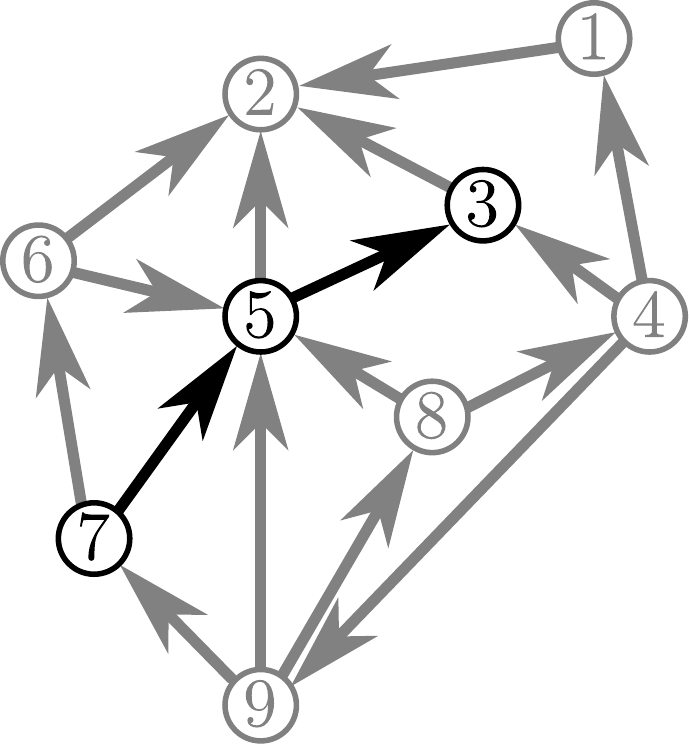}\label{fig:directconnect}}
\hspace{0.5cm}
\subfloat{{\footnotesize(b)}\hspace{-0.2cm}\includegraphics[width=3.8cm]{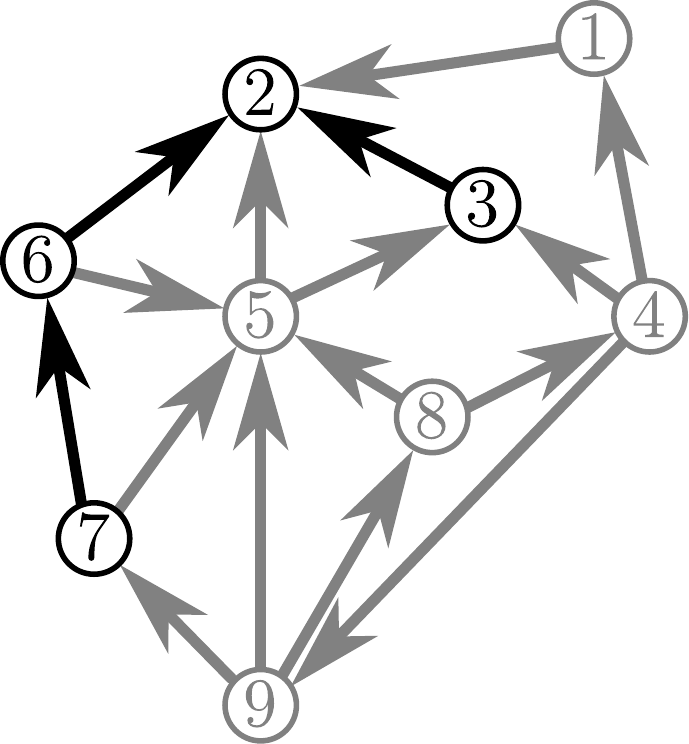}\label{fig:indirectconnect}}
\caption{A directed graph on $9$ nodes with: (a) a direct connection (i.e. a path) between nodes $3$ and $7$ highlighted, and (b) an indirect connection between nodes $3$ and $7$ highlighted.}
\label{fig:connections}
\end{figure}

We define a \emph{connection} in $\mathcal{G}$ between nodes $k$ and $j$ to consist of two paths, one starting at $k$ and the other at $j$ and which both terminate at the same node. A \emph{direct connection} between nodes $k$ and $j$ is a connection in which one path is trivial (i.e. either only node $k$ or only node $j$) - thus a direct connection is equivalent to a path. Conversely, an \emph{indirect connection} is one in which the terminal node of the two paths is neither node $k$ nor node $j$. Examples of direct and indirect connections are shown in Fig.~\ref{fig:connections}. A \emph{simple connection} is a connection that consists of two simple paths. Note that in both a connection and a simple connection, multiple nodes may be common between the two paths.

The graph $\mathcal{G}$ is \emph{connected} if it contains a globally reachable node $k$; i.e. there exists a node $k$ such that there is a path in $\mathcal{G}$ from $i$ to $k$ for every node $i$. Equivalently, $\mathcal{G}$ is connected if and only if a connection exists between any pair of nodes. The eigenvalue $0$ is simple if and only if $\cal G$ is connected \cite{Agaev05}. An undirected graph is connected if and only if there is a path between any pair of nodes.

A directed (respectively, undirected) \emph{tree} is a connected graph on $N$ nodes that contains exactly $N - 1$ directed (respectively, undirected) edges. A \emph{leaf} in a directed tree is any node with zero in-degree, and a leaf in an undirected tree is any node with only one neighbor. The \emph{root} of a directed tree is a node with zero out-degree; note that every directed tree will contain precisely one such node. A \emph{branch} of a directed tree is a path from a leaf to the root. 

Let $\Pi \mathrel{\mathop :}= I_N - \frac{1}{N}\mathbf{1}_N\mathbf{1}_N^T$ denote the orthogonal projection matrix onto the subspace of $\mathbb{R}^N$ perpendicular to $\mathbf{1}_N$. We will make a slight abuse of notation and use $\mathbf{1}_N^\perp$ to denote this subspace, instead of $\text{span}\!\left\{\mathbf{1}_N\right\}^\perp$. The matrix $\Pi$ is symmetric and since $L\mathbf{1}_N = \mathbf{0}$, $L\Pi = L$ and $\Pi L^T = L^T$ for any graph. Furthermore, $\Pi L = L$ for any balanced graph (including every undirected graph).

Let $Q \in \mathbb{R}^{(N-1)\times N}$ be a matrix whose rows form an orthonormal basis for $\mathbf{1}_N^\perp$. This is equivalent to requiring that $\begin{bmatrix}\frac{1}{\sqrt{N}}\mathbf{1}_N & Q^T\end{bmatrix}$ is an orthogonal matrix, or more explicitly,
\begin{equation}\label{eqn:propq}
Q\mathbf{1}_N = \mathbf{0}, \; QQ^T = I_{N-1} \text{ and } Q^TQ = \Pi.
\end{equation}
Using these properties, it follows that $Q\Pi = Q$ and $\Pi Q^T = Q^T$.

\section{An extended definition of effective resistance}\label{sec:def}
We now proceed to examine the derivation of effective resistance for undirected graphs and compare the matrices involved to those that arise in control-theoretic applications. Using this comparison, we propose a generalization of effective resistance to directed graphs that preserves key control-theoretic properties related to consensus-type dynamics.

A complete derivation of the standard notion of effective resistance is given in \cite{Klein93}, in which the effective resistance between two nodes in an undirected graph can be calculated by appropriately applying Kirchhoff's voltage and current laws. This calculation relies on what the authors call the ``generalized inverse'' of the Laplacian matrix, a matrix $X$ which satisfies
\begin{equation}\label{eqn:geninv}
XL = LX = \Pi \text{ and } X\Pi = \Pi X = X.
\end{equation}
Then, if we let $X = [x_{i,j}]$, the effective resistance is given by
\begin{equation}\label{eqn:rdef}
r_{k,j} = \left(\mathbf{e}_N^{(k)} \!-\! \mathbf{e}_N^{(j)}\right)^{\!\!T} \!\!X\! \left(\mathbf{e}_N^{(k)} \!-\! \mathbf{e}_N^{(j)}\right) = x_{k,k} \!+\! x_{j,j} \!-\! 2x_{k,j}.
\end{equation}
Although \refe{eqn:geninv} are not the defining equations for the Moore-Penrose generalized inverse, it is easy to show that for a symmetric Laplacian matrix $L$, any solution to \refe{eqn:geninv} will indeed be the Moore-Penrose generalized inverse of $L$ (as well as the group inverse of $L$). In fact, it is standard practice to define the effective resistance in terms of the Moore-Penrose generalized inverse \cite{Xiao03}. 

In \cite{Klein93}, the authors describe $X$ in the following way (with notation changed to match this paper):
\begin{definition}[Klein and Randi\'{c}, \cite{Klein93}]\label{def:undirX}
$X$ is equal on $\mathbf{1}_N^\perp$ to the inverse of $L$ and is otherwise $\Mzero$.
\end{definition}
It is therefore instructive to characterize the action of $L$ restricted to the subspace $\mathbf{1}_N^{\perp}$. Suppose we choose an orthonormal basis for $\mathbf{1}_N^\perp$ and let $Q \in \mathbb{R}^{\left(N-1\right)\times N}$ be the matrix formed with these basis vectors as rows. Then for any $\mathbf{v} \in \mathbb{R}^N$, $\overline{\mathbf{v}} \mathrel{\mathop :}= Q\mathbf{v}$ is a coordinate vector (with respect to the chosen basis) of the orthogonal projection of $\mathbf{v}$ onto $\mathbf{1}_N^\perp$ and for any $M \in \mathbb{R}^{N\times N}$, $\overline{M} \mathrel{\mathop :}= QMQ^T$ is the $\left(N-1\right)\times\left(N-1\right)$ matrix whose action on $\mathbf{1}_N^\perp$ is identical to $M$, in the sense that $\overline{M}\overline{\mathbf{v}} = \overline{M \mathbf{v}}$ for any $\mathbf{v} \in \mathbf{1}_N^\perp$.

Thus on $\mathbf{1}_N^\perp$, the Laplacian matrix is equivalent to
\begin{equation}\label{eqn:lbar}
\overline{L} = QLQ^T,
\end{equation}
which we refer to as the \emph{reduced Laplacian}. We can see that $\overline{L}$ is a symmetric matrix if and only if the graph is undirected. The reduced Laplacian has the same eigenvalues as $L$ except for a single $0$ eigenvalue \cite{Young10}. Hence, for a connected graph, $\overline{L}$ is invertible. For undirected graphs, this allows us to give an explicit construction for $X$ as
\begin{equation}\label{eqn:xinv}
X = Q^T\overline{L}^{-1}Q,
\end{equation}
which satisfies Definition \ref{def:undirX} since $X\mathbf{1}_N = \mathbf{0}$ and $\overline{X \mathbf{v}} = \overline{L}^{-1}\overline{\mathbf{v}}$ for any $\mathbf{v} \in \mathbb{R}^N$. Furthermore, we can use \refe{eqn:propq} and the fact that $L = L\Pi = \Pi L$ for undirected graphs to show that \refe{eqn:xinv} satisfies \refe{eqn:geninv} when the graph is undirected.

It should be noted that $\overline{L}$ is not unique, since it depends on the choice of $Q$. However, if $Q$ and $Q^\prime$ both satisfy \refe{eqn:propq}, we can define $P \mathrel{\mathop :}= Q^\prime Q^T$. Then $Q^\prime = PQ$ and $P$ is orthogonal. Hence $X^\prime \mathrel{\mathop :}= Q^{\prime T}\left(Q^\prime L Q^{\prime T}\right)^{-1} Q^\prime = Q^T P^T\left(P QLQ^T P^T\right)^{-1} P Q = X$ and thus the computation of $X$ in \refe{eqn:xinv} is independent of the choice of $Q$.

These multiple ways (the Moore-Penrose generalized inverse, Definition \ref{def:undirX}, \refe{eqn:geninv} and \refe{eqn:xinv}) to describe the matrix $X$ no longer agree when the graph is directed. While \refe{eqn:xinv} still satisfies Definition \ref{def:undirX}, it does not satisfy \refe{eqn:geninv} (specifically, $LX$ no longer equals $\Pi$). Furthermore, the Moore-Penrose generalized inverse satisfies neither \refe{eqn:geninv} nor \refe{eqn:xinv} for non-symmetric Laplacian matrices. Thus, instead of seeking to extend the notion of effective resistance to directed graphs using one of the above descriptions (which all arose through an analysis of electrical networks that were, by definition, undirected), we draw inspiration from a different context not as fundamentally tied to electrical networks.

In previous work on distributed consensus-formation and evidence-accumulation for decision-making (\cite{Young10,Young11,Poulakakis12}), effective resistances arose due to a correspondence between covariance matrices and the matrix $X$ as described above (in the case of undirected graphs). These applications involved stochastic systems evolving on graphs with dynamics driven by the Laplacian matrix, and covariance matrices were sought to describe the distribution of node values. For general (i.e. directed or undirected) graphs, these covariance matrices were computed using integrals of the form
\begin{equation}\label{eqn:lbarcov}
\Sigma_1 = \int_0^\infty{Q^T\Exp{-\overline{L} t}\Exp{-\overline{L}^T t}Q\;dt} \text{ and}
\end{equation}
\begin{equation}\label{eqn:lbarint}
\Sigma = \int_0^\infty{\Exp{-\overline{L} t}\Exp{-\overline{L}^T t}\;dt}.
\end{equation}
Now, we can observe that $\Sigma_1 = Q^T\Sigma Q$, and that $\Sigma$ can also be expressed as the solution to the Lyapunov equation \cite{Dullerud00}
\begin{equation}\label{eqn:lyap}
\overline{L}\Sigma + \Sigma \overline{L}^T = I_{N-1}.
\end{equation}
It should be noted that \refe{eqn:lyap} has a unique positive definite solution when all the eigenvalues of $\overline{L}$ have positive real part (i.e. when the graph is connected) \cite{Dullerud00}. It is then clear that for undirected graphs (where $\overline{L}$ is symmetric),
\begin{equation}\label{eqn:siginv}
\Sigma = \frac{1}{2}\overline{L}^{-1},
\end{equation}
and so (using \refe{eqn:xinv}),
\[
\Sigma_1 = \frac{1}{2} X.
\]

It is this relationship that links these covariance matrices to the generalized inverse $X$, and hence to effective resistances. Since these covariance matrices arise naturally from directed graphs as well as undirected graphs, we use their solutions to define effective resistances on directed graphs. Thus, for any connected digraph, we let $\Sigma$ be the unique solution to the Lyapunov equation \refe{eqn:lyap}. Then, we let
\begin{equation}\label{eqn:xdef}
X \mathrel{\mathop :}= 2Q^T \Sigma Q,
\end{equation}
and notice that $X$ will be symmetric for any graph because $\Sigma$ is always symmetric. Finally we can use \refe{eqn:rdef} to define the effective resistance between any two nodes in the graph.

\begin{definition}\label{def:dirres}
Let $\mathcal{G}$ be a connected graph with $N$ nodes and Laplacian matrix $L$. Then the \emph{effective resistance between nodes $k$ and $j$ in $\mathcal{G}$} is defined as
\begin{align}\label{eqn:dirres}
r_{k,j}	&= \left(\mathbf{e}_N^{(k)} - \mathbf{e}_N^{(j)}\right)^T X \left(\mathbf{e}_N^{(k)} - \mathbf{e}_N^{(j)}\right) \nonumber\\
		&= x_{k,k} + x_{j,j} - 2x_{k,j},
\end{align}
where 
\begin{equation}\label{eqn:resdetails}
\begin{gathered}
X = 2Q^T\Sigma Q,\\
\overline{L}\Sigma + \Sigma\overline{L}^T = I_{N-1}, \\
\overline{L} = QLQ^T,
\end{gathered}
\end{equation}
and $Q$ is a matrix satisfying \refe{eqn:propq}.
\end{definition}

By summing all distinct effective resistances in a graph, we define the \emph{Kirchhoff index}, $K_f$, of the graph, that is,
\begin{equation}\label{eqn:kirchhoff}
K_f \mathrel{\mathop :}= \sum_{k < j}{r_{k,j}}.
\end{equation}
Our definition of $K_f$  generalizes the Kirchhoff index defined for undirected graphs \cite{Xiao03}.

\begin{remark}
In previous work (\cite{Young10, Young11}) we noted that an $\mathcal{H}_2$ norm $H$ related to linear consensus can be expressed in terms of the Kirchhoff index (for an undirected graph) as
\begin{equation}\label{eqn:H2}
H = \left(\frac{K_f}{2N}\right)^\frac{1}{2},
\end{equation}
since this $\mathcal{H}_2$ norm is derived from the trace of the matrix $\Sigma$ which solves \refe{eqn:lyap}. Therefore, if we take \refe{eqn:kirchhoff} to be the definition of the Kirchhoff index for directed graphs as well, we find that \refe{eqn:H2} continues to hold. Thus our definition of effective resistance immediately connects to the robustness to noise of linear consensus on directed graphs. In a similar fashion, the variance of each node in a balanced network of stochastic decision-makers can be computed (in part) using the diagonal entries of the matrix $\Sigma_1$ from \refe{eqn:lbarcov}. Then the variance of any particular node can be found using this definition of effective resistance \cite{Poulakakis12}.
\end{remark}

\section{Basic properties of our definition}\label{sec:prop}
Although Definition \ref{def:dirres} ensures that effective resistance maintains our desired relationship with some control-theoretic properties, such as those discussed above, by itself it remains an algebraic construction that yields little insight into the ways in which effective resistance depends on the graph structure. We now proceed to analyze our definition to understand some of its fundamental properties. In Section \ref{subsec:welldef} we verify that Definition \ref{def:dirres} results in a well-defined property of pairs of nodes in a connected digraph. In Section \ref{subsec:connect} we investigate how effective resistances depend on connections in the graph and extend Definition \ref{def:dirres} further to apply to disconnected graphs. Finally in Section \ref{subsec:metric} we determine that effective resistance is a distance-like function and explore the limitations of the triangle inequality for effective resistances in directed graphs.

\subsection{Effective resistance is well-defined}\label{subsec:welldef}
By construction, \refe{eqn:dirres} will yield the regular effective resistance for any undirected graph. However, we must confirm that our concept of effective resistance for directed graphs is well-defined. This is achieved by the following two lemmas.

\begin{lemma}\label{lem:indofq}
The value of the effective resistance between two nodes in a connected digraph is independent of the choice of $Q$.
\end{lemma}

\begin{IEEEproof}
Let $Q$ and $Q^\prime$ be two matrices that satisfy \refe{eqn:propq}, and let $r_{k,j}$ and $r_{k,j}^\prime$ be the corresponding effective resistances between nodes $k$ and $j$, computed by using $Q$ and $Q^\prime$ in \refe{eqn:dirres}, respectively. Furthermore, let $W \mathrel{\mathop :}= Q^\prime Q^T$. Then by \refe{eqn:propq}, $Q^\prime = W Q$ and $W W^T = W^T W = I_{N-1}$. Now, we can use \refe{eqn:lbar} and the properties of $W$ to write $\overline{L}^\prime = W\overline{L}W^T$.

Next, substituting this expression into \refe{eqn:lyap} for $\overline{L}^\prime$, we see that $\Sigma^\prime = W\Sigma W^T$. Finally, we find that $X^\prime = 2Q^{\prime T}\Sigma^\prime Q^\prime = X,$ and hence, $r_{k,j}^\prime = r_{k,j}$.
\end{IEEEproof}

From Lemma \ref{lem:indofq} we see that the no matter how it is computed, the effective resistance between two nodes will be the same unique number. Next, we will show in Lemma \ref{lem:relabel} that the effective resistance is a property of a pair of nodes, irrespective of the way in which they are labelled.

\begin{lemma}\label{lem:relabel}
The value of the effective resistance between two nodes in a connected digraph is independent of the labeling of the nodes.
\end{lemma}

\begin{IEEEproof}
Any two labelings of the nodes in a graph can be related via a permutation. Suppose $L$ and $L^\prime$ are two Laplacian matrices associated with the same graph, but with different labelings of the nodes. Then $L^\prime$ can be found from $L$ by permuting its rows and columns. That is, there exists an $N\times N$ permutation matrix $P$ such that $L^\prime = P L P^T$. Note that as a permutation matrix, there is exactly one $1$ in every row and column of $P$ with every other entry equal to $0$. Furthermore, $P^{-1} = P^T$, $P\mathbf{1}_N = \mathbf{1}_N$ and $\mathbf{1}_N^T P = \mathbf{1}_N^T$. Thus we can observe that $QP = \overline{P}Q,$ $PQ^T = Q^T\overline{P}$ and $\overline{P}^{-1} = \overline{P}^T,$ where, as usual, $\overline{P} = QPQ^T$.

Now, we can use \refe{eqn:lbar} and the properties of $P$ to write $\overline{L}^\prime = \overline{P}\,\overline{L}\,\overline{P}^T$. Then the solution to the Lyapunov equation associated with $\overline{L}^\prime$ becomes $\Sigma^\prime = \overline{P} \Sigma \overline{P}^T$. Hence, we observe that $X^\prime	= P X P^T$.

Thus if $P$ permutes node $k$ to node $m$ and node $j$ to node $\ell$, we find that $r^\prime_{m,\ell} = r_{k,j}$.
\end{IEEEproof}

\subsection{Effective resistance depends on connections between nodes}\label{subsec:connect}
\begin{figure}
\centering
\includegraphics[width=3.8cm]{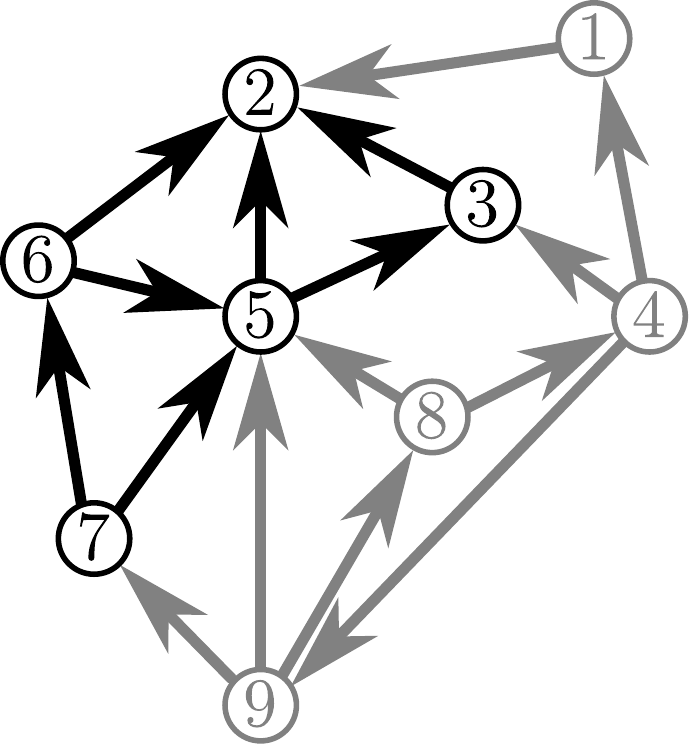}
\caption{A directed graph on $9$ nodes with the reachable subgraph of nodes $3$ and $7$ highlighted.}
\label{fig:reachable}
\end{figure}

Next we consider which features of a digraph will affect the effective resistance between a given pair of nodes. For undirected graphs, we know that effective resistances depend on every path between a pair of nodes \cite{Klein93}. The situation becomes more complicated with directed graphs since there can exist pairs of nodes in a connected digraph with no path between them. Instead of looking at paths between nodes, we therefore have to consider connections. To incorporate all of the connections between two nodes, we examine the \emph{reachable subgraph}, an example of which is shown in Fig.~\ref{fig:reachable}.

\begin{definition}\label{def:reachable}
The \emph{reachable subgraph, denoted $\mathcal{R}_\mathcal{G}(k,j)$, of nodes $k$ and $j$ in the graph $\mathcal{G}$} is the graph formed by every node in $\mathcal{G}$ that is reachable from node $k$ or node $j$ and every edge in $\mathcal{G}$ between these nodes.
\end{definition}

As we demonstrate in the following lemma, if $\mathcal{G}$ is connected the reachable subgraph of nodes $k$ and $j$ is precisely the subgraph formed by every connection between them.

\begin{lemma}\label{lem:reachable}
If $\mathcal{G}$ is connected, then for any pair of nodes $k$ and $j$,
\begin{enumerate}[(i) ]
\item $\mathcal{R}_\mathcal{G}(k,j)$ is connected,
\item Every node in $\mathcal{R}_\mathcal{G}(k,j)$ is part of a connection between nodes $k$ and $j$,
\item Every edge in $\mathcal{R}_\mathcal{G}(k,j)$ is part of a connection between nodes $k$ and $j$ and
\item Every connection between nodes $k$ and $j$ is contained in $\mathcal{R}_\mathcal{G}(k,j)$.
\end{enumerate}
\end{lemma}

\begin{IEEEproof}
\begin{enumerate}[(i) ]
\item Since $\mathcal{G}$ is connected, there is a node, $\ell$, in $\mathcal{G}$ which is reachable from every other node. Since $\ell$ is reachable from nodes $k$ and $j$, it is also in $\mathcal{R}_\mathcal{G}(k,j)$. Now, suppose that $m$ is a node in $\mathcal{R}_\mathcal{G}(k,j)$. Then there is a path in $\mathcal{G}$ from $m$ to $\ell$. Since $m$ is reachable from either $k$ or $j$, every node along this path is as well. Thus, this path is contained in $\mathcal{R}_\mathcal{G}(k,j)$ and so $\ell$ is reachable (in $\mathcal{R}_\mathcal{G}(k,j)$) from every node in $\mathcal{R}_\mathcal{G}(k,j)$.

\item Let $m$ be a node in $\mathcal{R}_\mathcal{G}(k,j)$. Then $m$ must be reachable from either $k$ or $j$. Without loss of generality, suppose that $m$ is reachable from $k$. Then (as we saw in part (i)) there must be a path in $\mathcal{R}_\mathcal{G}(k,j)$ from $m$ to the globally reachable node $\ell$ as well as a path from $j$ to $\ell$. Thus $m$ is part of a connection between $k$ and $j$.

\item Let $\left(m,n\right)$ be an edge in $\mathcal{R}_\mathcal{G}(k,j)$. Without loss of generality, suppose that $m$ is reachable from $k$. Then $n$ is also reachable from $k$. Then (as we saw in part (i)) there must be a path in $\mathcal{R}_\mathcal{G}(k,j)$ from $n$ to the globally reachable node $\ell$ as well as a path from $j$ to $\ell$. Thus $\left(m,n\right)$ is part of a connection between $k$ and $j$.

\item Every node along a path is reachable from the node where the path started. Thus, every node in a connection between $k$ and $j$ is reachable from either $k$ or $j$ and hence in $\mathcal{R}_\mathcal{G}(k,j)$. Since $\mathcal{R}_\mathcal{G}(k,j)$ contains every edge in $\mathcal{G}$ between its nodes, every edge in the connection must also be in $\mathcal{R}_\mathcal{G}(k,j)$.
\end{enumerate}
\end{IEEEproof}

Next we proceed to show in Theorem \ref{theo:rconnect} that the effective resistance between two nodes in a connected digraph can only depend on the connections between them. The proof relies on the following lemma, which describes sufficient conditions under which effective resistances in a subgraph will be equal to those in the original graph.

\begin{lemma}\label{lem:equalX}
Suppose that $\mathcal{G}_1$ is a connected subgraph (containing $N_1$ nodes and with Laplacian matrix $L_1$) of a connected graph $\mathcal{G}$ (containing $N$ nodes and with Laplacian matrix $L$) and suppose that the nodes in $\mathcal{G}_1$ are labelled $1$ through $N_1$. Let $Q_1 \in \mathbb{R}^{(N_1-1)\times N_1}$ be a matrix satisfying \refe{eqn:propq} and suppose there is a $Q \in \mathbb{R}^{(N-1)\times N}$ satisfying \refe{eqn:propq} that can be written as
\[
Q = \begin{bmatrix}Q_1 & \Mzero\\\alpha\mathbf{1}_{N_1}^T & \mathbf{r}^T\\\Mzero & S\end{bmatrix}
\]
for some $\alpha \in \mathbb{R}$, $\mathbf{r} \in \mathbb{R}^{N-N_1}$ and $S \in \mathbb{R}^{(N-N_1 - 1)\times(N-N_1)}$. If the solution to \refe{eqn:lyap} for $\mathcal{G}$ (with $\overline{L} = QLQ^T$) can be written as
\[
\Sigma = \begin{bmatrix}\Sigma_1 & \mathbf{t} & U\\\mathbf{t}^T & v & \mathbf{w}^T\\U^T & \mathbf{w} & Y\end{bmatrix},
\]
for some $\mathbf{t} \in \mathbb{R}^{N_1-1}$, $U \in \mathbb{R}^{(N_1-1)\times(N-N_1-1)}$, $v \in \mathbb{R}$, $\mathbf{w} \in \mathbb{R}^{N-N_1-1}$ and $Y \!\in \mathbb{R}^{(N-N_1-1)\times(N-N_1-1)}$ with $Y = Y^T$ and where $\Sigma_1 \in \mathbb{R}^{(N_1-1)\times(N_1-1)}$ is the solution to \refe{eqn:lyap} for $\mathcal{G}_1$ (with $\overline{L}_1 = Q_1L_1Q_1^T$), then for any $k, j \leq N_1$, the effective resistance between nodes $k$ and $j$ in $\mathcal{G}$ is equal to the effective resistance between the same two nodes in $\mathcal{G}_1$.
\end{lemma}

\begin{IEEEproof}
Effective resistances in $\mathcal{G}_1$ can be found from $X_1 = 2Q_1^T\Sigma_1 Q_1$ as
\[
r_{1\,k,j} = x_{1\,k,k} + x_{1\,j,j} - 2x_{1\,k,j}.
\]
To compute effective resistances in $\mathcal{G}$, we must examine $X = 2Q^T\Sigma Q$. Using the matrices given above, we obtain
\begin{multline*}
X = \left[\begin{matrix}2Q_1^T\Sigma_1Q_1 \!+\! 2\alpha\mathbf{1}_{N_1}\mathbf{t}^TQ_1 \!+\! 2\alpha Q_1^T\mathbf{t}\mathbf{1}_{N_1}^T \!+\! 2\alpha^2 v \mathbf{1}_{N_1}\mathbf{1}_{N_1}^T \\ \ast \end{matrix}\right.\\
\left.\begin{matrix} \ast \\ \ast \end{matrix}\right].
\end{multline*}
If we let $\mathbf{p} \mathrel{\mathop :}= 2\alpha Q_1^T\mathbf{t} = [p_i]$, we can write
\[
X = \begin{bmatrix}X_1 + \mathbf{1}_{N_1}\mathbf{p}^T + \mathbf{p}\mathbf{1}_{N_1}^T + 2\alpha^2 v \mathbf{1}_{N_1}\mathbf{1}_{N_1}^T & \ast \\ \ast & \ast \end{bmatrix}.
\]

Finally, since nodes $k$ and $j$ are both in $\mathcal{G}_1$, we obtain
\begin{align*}
r_{k,j}	&= x_{1\,k,k} + 2p_k + 2\alpha^2 v + x_{1\,j,j} + 2p_j + 2\alpha^2 v - 2x_{1\,k,j} \\
		&\hspace{5.0cm} {} - 2(p_k+p_j) - 4\alpha^2 v \\
		&= r_{1\,k,j}.
\end{align*}

Note that the same calculation applies if $N = N_1 + 1$, in which case the $\Mzero_{(N-N_1-1)\times N_1}$, $S$, $U$, $\mathbf{w}$ and $Y$ blocks of $Q$ and $\Sigma$ are all empty.
\end{IEEEproof}

Now we can state our first main result.

\begin{theorem}\label{theo:rconnect}
The effective resistance between nodes $k$ and $j$ in a connected graph $\mathcal{G}$ is equal to the effective resistance between nodes $k$ and $j$ in $\mathcal{R}_\mathcal{G}(k,j)$.
\end{theorem}

\begin{IEEEproof}
Let $\mathcal{G}_1 = \mathcal{R}_\mathcal{G}(k,j)$. Let $N_1$, $A_1$, $D_1$ and $L_1$ be the number of nodes, the adjacency matrix, the matrix of node out-degrees and the Laplacian matrix of $\mathcal{G}_1$, respectively. Let $Q_1 \in \mathbb{R}^{(N_1-1)\times N_1}$ satisfy \refe{eqn:propq}. Since $\mathcal{G}_1$ is connected by Lemma \ref{lem:reachable}, we can find matrices $\overline{L}_1$, $\Sigma_1$ and $X_1$ from \refe{eqn:resdetails} for $\mathcal{G}_1$. 

Let $\mathcal{G}_2$ be the subgraph of $\mathcal{G}$ formed by every node in $\mathcal{G}$ which is not in $\mathcal{G}_1$ and every edge in $\mathcal{G}$ between these nodes. Then $\mathcal{G}_2$ will contain $N_2$ nodes and have associated matrices $A_2$, $D_2$, $L_2$, $Q_2$ and $\overline{L}_2$.

Now, if there was an edge $\left(m,n\right)$ in $\mathcal{G}$ from a node $m$ in $\mathcal{G}_1$ to a node $n$ in $\mathcal{G}_2$, then $n$ would be reachable from either $k$ or $j$ (as $m$ is reachable from one of these nodes). Thus there are no edges in $\mathcal{G}$ leading from nodes in $\mathcal{G}_1$ to nodes in $\mathcal{G}_2$. Therefore, if we order the nodes in $\mathcal{G}$ with the nodes in $\mathcal{G}_1$ listed first, followed by the nodes in $\mathcal{G}_2$, then the adjacency matrix of $\mathcal{G}$ can be written as
\[
A = \begin{bmatrix}A_1 & \Mzero\\A_{2,1} & A_2\end{bmatrix},
\]
where $A_{2,1} \in \mathbb{R}^{N_2 \times N_1}$ contains the edge weights for all edges in $\mathcal{G}$ leading from nodes in $\mathcal{G}_2$ to nodes in $\mathcal{G}_1$. Similarly, we can write the matrix of node out-degrees as
\[
D = \begin{bmatrix}D_1 & \Mzero\\\Mzero & D_2 + D_{2,1}\end{bmatrix},
\]
where $D_{2,1}$ is the diagonal matrix containing the row sums of $A_{2,1}$ along its diagonal. Utilizing these two expressions, we find that the Laplacian matrix of $\mathcal{G}$ can be written as
\[
L = \begin{bmatrix}L_1 & \Mzero\\-A_{2,1} & L_2+D_{2,1}\end{bmatrix}.
\]

Now, let
\[
Q = \begin{bmatrix}Q_1 & \Mzero\\ \alpha \mathbf{1}_{N_1}^T & -\beta\mathbf{1}_{N_2}^T \\\Mzero & Q_2\end{bmatrix},
\]
where $\alpha = \sqrt{\dfrac{N_2}{N_1\left(N_1+N_2\right)}}$ and $\beta = \sqrt{\dfrac{N_1}{N_2\left(N_1+N_2\right)}}$. Then $Q$ satisfies \refe{eqn:propq} (note that $N = N_1 + N_2$). Substituting this matrix into \refe{eqn:lbar} gives us
\begin{multline*}
\overline{L} = \left[\begin{matrix}\overline{L}_1 & \mathbf{0} \\ \alpha\mathbf{1}_{N_1}^T L_1 Q_1^T + \beta \mathbf{1}_{N_2}^T A_{2,1} Q_1^T & \beta\left(\alpha + \beta\right)\mathbf{1}_{N_2}^T \mathbf{d}_{2,1} \\-Q_2 A_{2,1} Q_1^T & -\left(\alpha + \beta\right)Q_2 \mathbf{d}_{2,1} \end{matrix}\right.\\
\left.\begin{matrix} \Mzero\\  -\beta \mathbf{1}_{N_2}^T L_2 Q_2^T - \beta \mathbf{1}_{N_2}^T D_{2,1} Q_2^T\\ \overline{L}_2 + Q_2 D_{2,1} Q_2^T\end{matrix}\right],
\end{multline*}
where $\mathbf{d}_{2,1} \mathrel{\mathop :}= D_{2,1} \mathbf{1}_{N_2} = A_{2,1} \mathbf{1}_{N_1}$.

In order to compute effective resistances in $\mathcal{G}$, we must find the matrix $\Sigma$ which solves \refe{eqn:lyap}. Since we have partitioned $\overline{L}$ into a $3\times 3$ block matrix, we will do the same for $\Sigma$. Let
\[
\Sigma = \begin{bmatrix} S & \mathbf{t} & U\\ \mathbf{t}^T & v & \mathbf{w}^T\\U^T & \mathbf{w} & Y\end{bmatrix},
\]
where $S \in \mathbb{R}^{(N_1 - 1)\times(N_1 - 1)}$ and $Y \in \mathbb{R}^{(N_2 - 1)\times(N_2 - 1)}$ are symmetric matrices, $U \in \mathbb{R}^{(N_1 - 1)\times(N_2 - 1)}$, $\mathbf{t} \in \mathbb{R}^{N_1 - 1}$, $\mathbf{w} \in \mathbb{R}^{N_2-1}$ and $v \in \mathbb{R}$. Then when we multiply out the matrices of \refe{eqn:lyap} and equate the $\left(1,1\right)$ blocks, we find
\[
\overline{L}_1 S + S\overline{L}_1^T = I_{N_1-1},
\]
which implies that
\[
S = \Sigma_1.
\]

Thus, by Lemma \ref{lem:equalX}, the effective resistance between two nodes in $\mathcal{G}_1$ is equal to the effective resistance between the same two nodes in $\mathcal{G}$.
\end{IEEEproof}

We can use Theorem \ref{theo:rconnect} to partially extend the definition of effective resistance to disconnected digraphs. To do this, we will first define \emph{connection subgraphs}.

\begin{definition}\label{def:connection}
A \emph{connection subgraph between nodes $k$ and $j$ in the graph $\mathcal{G}$} is a maximal connected subgraph of $\mathcal{G}$ in which every node and edge form part of a connection between nodes $k$ and $j$ in $\mathcal{G}$. That is, a connection subgraph is formed from the union of connections between nodes $k$ and $j$, and the addition of any other connections would make the subgraph disconnected. If only one connection subgraph exists in $\mathcal{G}$ between nodes $k$ and $j$, it is referred to as \emph{the} connection subgraph and is denoted by $\mathcal{C}_\mathcal{G}(k,j)$.
\end{definition}

\begin{figure}
\centering
\subfloat{{\footnotesize(a)}\hspace{0.2cm}\includegraphics[width=1.5cm]{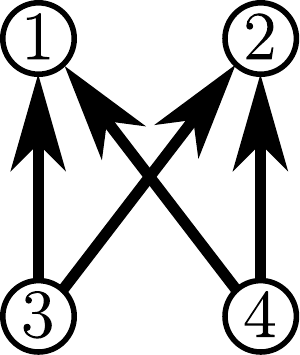}\label{fig:disconnectgraph}}
\hspace{0.3cm}
\subfloat{{\footnotesize(b)}\hspace{0.2cm}\includegraphics[width=1.5cm]{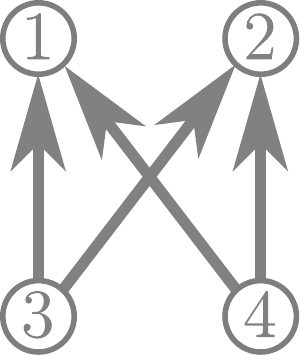}\label{fig:disconnect12}}
\hspace{0.3cm}
\subfloat{{\footnotesize(c)}\hspace{0.2cm}\includegraphics[width=1.5cm]{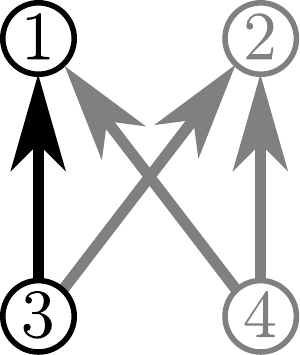}\label{fig:disconnect13}}
\hspace{0.3cm}
\subfloat{{\footnotesize(d)}\hspace{0.2cm}\includegraphics[width=1.5cm]{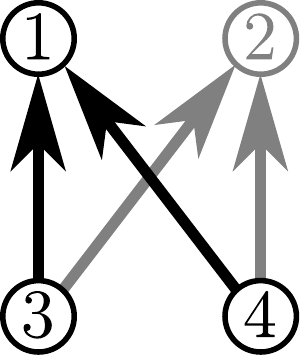}\label{fig:disconnect341}}
\hspace{0.3cm}
\subfloat{{\footnotesize(e)}\hspace{0.2cm}\includegraphics[width=1.5cm]{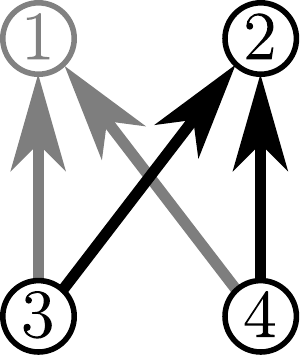}\label{fig:disconnect342}}
\caption{(a) A disconnected graph on $4$ nodes. (b) There are no connection subgraphs between nodes $1$ and $2$. (c) The connection subgraph between nodes $1$ and $3$ is highlighted. (d) One connection subgraph between nodes $3$ and $4$ is highlighted. (e) A second connection subgraph between nodes $3$ and $4$ is highlighted. In this example, the effective resistance between nodes $3$ and $4$ is undefined.}
\label{fig:disconnected}
\end{figure}

From Lemma \ref{lem:reachable} we know that $\mathcal{C}_{\mathcal{G}}(k,j) = \mathcal{R}_{\mathcal{G}}(k,j)$ if $\mathcal{G}$ is connected. However, a disconnected graph may contain $0$, $1$ or more connection subgraphs between a pair of nodes. There will be no connection subgraphs precisely when there are no connections between nodes $k$ and $j$ in $\mathcal{G}$. However, there may also be multiple connections between a pair of nodes that lead to multiple connection subgraphs. A simple example of connection subgraphs in a disconnected graph is shown in Fig.~\ref{fig:disconnected}

By definition, whenever it exists, $\mathcal{C}_{\mathcal{G}}(k,j)$ will be connected, and we can thus compute effective resistances in it. We can now define effective resistances between some node pairs in any digraph, whether or not it is connected, as follows

\begin{definition}\label{def:generalres}
The effective resistance between nodes $k$ and $j$ in a graph $\mathcal{G}$ is
\[
r_{k,j} = \begin{cases} \infty &\text{if there are no connections}\\ & \text{between nodes $k$ and $j$} \\ r_{k,j} \text{ in } \mathcal{C}_{\mathcal{G}}(k,j) & \text{if $\mathcal{C}_\mathcal{G}(k,j)$ exists}\\  &\text{(computed using \refe{eqn:dirres})}\\ \text{undefined} & \text{otherwise.}\end{cases}
\]
\end{definition}

By Theorem \ref{theo:rconnect}, this new definition specializes to our original definition of effective resistance for connected graphs. For certain applications, there may be an appropriate way to handle pairs of nodes with multiple connection subgraphs, but that falls outside the scope of the present work.

\begin{figure}
\centering
\subfloat{{\footnotesize(a)}\hspace{0.5cm}\includegraphics[width=0.6cm]{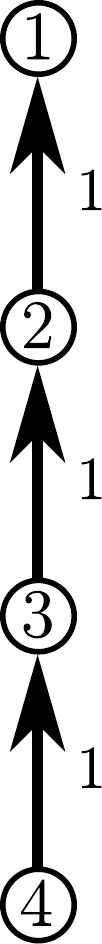}\label{fig:pathcounter}}
\hspace{2cm}
\subfloat{{\footnotesize(b)}\hspace{0.5cm}\includegraphics[width=0.6cm]{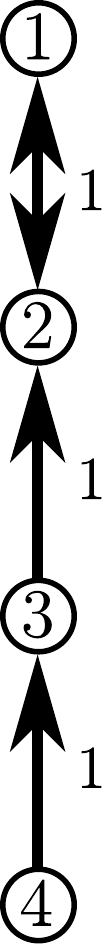}\label{fig:linecounter}}
\caption{Two simple $4$-node graphs: (a) $\mathcal{G}^\text{path}_4$, a $4$-node directed path graph with unit edge weights, and (b) $\mathcal{G}^\text{line}$, similar to $\mathcal{G}^\text{path}_4$ but with the directed edge $(2,1)$ replaced by an undirected edge.}
\label{fig:pathcounterexample}
\end{figure}

In undirected graphs, we know that the effective resistance between two nodes does not depend on edges that do not lie on any simple path between the nodes \cite{Klein93}. Unfortunately, the situation is not as straightforward for directed graphs. Consider the $4$-node graphs shown in Fig.~\ref{fig:pathcounterexample}. In $\mathcal{G}^\text{path}_4$, we can compute that $r_{3,4} = 2$. However in $\mathcal{G}^\text{line}$, even though no new simple connections were introduced between nodes $3$ and $4$, the effective resistance is now $r_{3,4} = \frac{16}{9}$.

Despite this (perhaps) unexpected behavior, we are able to show below that certain parts of the connection subgraph do not affect the effective resistance between two nodes. The proof relies on the following lemma, which provides a solution to Lyapunov equations with a certain structure.

\begin{lemma}\label{lem:lyapsol}
Suppose that $L_1 \in \mathbb{R}^{N_1 \times N_1}$ is the Laplacian matrix of a connected graph satisfying $L_1^T \mathbf{e}^{(1)}_{N_1} = \mathbf{0}$ and $Q_1 \in \mathbb{R}^{(N_1-1)\times N_1}$ satisfies \refe{eqn:propq}. Let $\alpha = \dfrac{1}{\sqrt{N_1(N_1+1)}}$, $\beta = \dfrac{N_1}{\sqrt{N_1(N_1+1)}}$ and $a > 0$. If $\overline{L}_1 = Q_1 L_1 Q_1^T$, then a solution to the Lyapunov equation
\begin{multline*}
\begin{bmatrix}\overline{L}_1 + aQ_1\mathbf{e}^{(1)}_{N_1}\mathbf{e}^{(1) T}_{N_1}Q_1^T & a(\alpha+\beta)Q_1\mathbf{e}_{N_1}^{(1)} \\ \alpha\mathbf{1}_{N_1}^T L_1 Q_1^T + a\alpha \mathbf{e}_{N_1}^{(1) T} Q_1^T & \frac{a}{N_1}\end{bmatrix} \begin{bmatrix} S & \mathbf{t}\\ \mathbf{t}^T & u\end{bmatrix} \\
{}\!+\! \begin{bmatrix} S & \mathbf{t}\\ \mathbf{t}^T & u\end{bmatrix}\!\! \begin{bmatrix}\overline{L}_1^T \!+\! aQ_1\mathbf{e}^{(1)}_{N_1}\mathbf{e}^{(1) T}_{N_1}Q_1^T & \alpha Q_1 L_1^T \mathbf{1}_{N_1} \!+\! a\alpha Q_1 \mathbf{e}_{N_1}^{(1)}\\ a(\alpha+\beta)\mathbf{e}_{N_1}^{(1) T}Q_1^T & \frac{a}{N_1}\end{bmatrix} \\
= \begin{bmatrix}I_{N_1-1} & \mathbf{0}\\\mathbf{0}^T & 1\end{bmatrix}
\end{multline*}
is
\begin{equation}\label{eqn:lyapsol}
\begin{gathered}
S = \Sigma_1, \\
\mathbf{t} = -N_1 \alpha \Sigma_1 Q_1 \mathbf{e}^{(1)}_{N_1} \text{ and} \\
u = \frac{N_1}{2a} + \frac{N_1^2 \alpha^2}{a}\left(\mathbf{1}_{N_1}^T L_1 + a\mathbf{e}^{(1) T}_{N_1}\right)Q_1^T\Sigma_1 Q_1 \mathbf{e}^{(1)}_{N_1},
\end{gathered}
\end{equation}
where $\Sigma_1$ is the solution to \refe{eqn:lyap} for $\overline{L}_1$.
\end{lemma}

\begin{IEEEproof}
First we note that $\Sigma_1$ exists since $L_1$ is the Laplacian of a connected graph. Next, equating blocks of the given matrix equation gives us
\begin{align}
\overline{L}_1 S + S \overline{L}_1^T + aQ_1\mathbf{e}^{(1)}_{N_1}\mathbf{e}^{(1) T}_{N_1}Q_1^T S + a S Q_1\mathbf{e}^{(1)}_{N_1}\mathbf{e}^{(1) T}_{N_1}Q_1^T \hspace{-1.6cm}&\nonumber\\
 {} + a(\alpha + \beta)\left(Q_1\mathbf{e}_{N_1}^{(1)}\mathbf{t}^T + \mathbf{t}\mathbf{e}_{N_1}^{(1) T}Q_1^T\right)	&= I_{N_1-1},\hspace{-0.0cm} \label{eqn:lyapsolS}\\
\alpha S Q_1 \left(L_1^T \mathbf{1}_{N_1} + a \mathbf{e}^{(1)}_{N_1}\right) + a u (\alpha + \beta)Q_1\mathbf{e}_{N_1}^{(1)}  \hspace{-0.3cm}& \nonumber\\
 {} + \left(\frac{a}{N_1}I_{N_1-1} + \overline{L}_1 + aQ_1\mathbf{e}^{(1)}_{N_1}\mathbf{e}^{(1) T}_{N_1}Q_1^T\right)\mathbf{t}	&= \mathbf{0} \text{ and} \label{eqn:lyapsolt}\\
\frac{2a u}{N_1} + 2\alpha\left(\mathbf{1}_{N_1}^T L_1 + a \mathbf{e}_{N_1}^{(1) T}\right)Q_1^T\mathbf{t}	&= 1. \label{eqn:lyapsolu}
\end{align}

By directly substituting \refe{eqn:lyapsol} into \refe{eqn:lyapsolS} and \refe{eqn:lyapsolu} (and noting that $N_1 \alpha(\alpha + \beta) = 1$), we observe that \refe{eqn:lyapsol} satisfies \refe{eqn:lyapsolS} and \refe{eqn:lyapsolu}. Hence, we now focus our attention on \refe{eqn:lyapsolt}. Substituting \refe{eqn:lyapsol} into the left-hand side of \refe{eqn:lyapsolt} gives us
\begin{multline}\label{eqn:lhs1}
LHS = \alpha \Sigma_1 Q_1 L_1^T \mathbf{1}_{N_1} - N_1 \alpha \overline{L}_1 \Sigma_1 Q_1 \mathbf{e}^{(1)}_{N_1} + \frac{1}{2\alpha}Q_1 \mathbf{e}^{(1)}_{N_1} \\
{} + N_1 \alpha Q_1 \mathbf{e}^{(1)}_{N_1}\mathbf{1}_{N_1}^T L_1 Q_1^T \Sigma_1 Q_1 \mathbf{e}^{(1)}_{N_1}.
\end{multline}

By \refe{eqn:lyap}, we know that $\overline{L}_1 \Sigma_1 = I_{N_1-1} - \Sigma_1 \overline{L}_1^T$. Therefore, using \refe{eqn:propq} and our assumption that $L_1^T \mathbf{e}^{(1)}_{N_1} = \mathbf{0}$, we have
\[
\overline{L}_1\Sigma_1 Q_1 \mathbf{e}^{(1)}_{N_1}	 = Q_1 \mathbf{e}^{(1)}_{N_1}	 + \frac{1}{N_1}\Sigma_1Q_1L_1^T\mathbf{1}_{N_1}.
\]
Thus $N_1 \alpha \overline{L}_1 \Sigma_1 Q_1 \mathbf{e}^{(1)}_{N_1} = N_1\alpha Q_1 \mathbf{e}^{(1)}_{N_1} + \alpha\Sigma_1Q_1L_1^T\mathbf{1}_{N_1}$, and \refe{eqn:lhs1} becomes
\begin{equation}\label{eqn:lhs2}
LHS = \left(\frac{1}{2\alpha} - N_1 \alpha + N_1 \alpha\mathbf{1}_{N_1}^T L_1 Q_1^T \Sigma_1 Q_1 \mathbf{e}^{(1)}_{N_1}\right)Q_1 \mathbf{e}^{(1)}_{N_1}.
\end{equation}

Next, if we define $V$ to be the matrix $V \mathrel{\mathop:}= L_1 Q_1^T \Sigma_1 Q_1 + Q_1^T \Sigma_1 Q_1 L_1^T$, we have that $V = V^T$ and (using \refe{eqn:propq}),
\[
\mathbf{1}_{N_1}^T V \mathbf{e}^{(1)}_{N_1} = \mathbf{1}_{N_1}^T L_1 Q_1^T \Sigma_1 Q_1 \mathbf{e}^{(1)}_{N_1}.
\]
But pre- and post-multiplying $V$ by $\Pi = Q_1^T Q_1$ and using \refe{eqn:lyap} gives us $\Pi V \Pi = \Pi$, and then by pre- and post-multiplying by $\mathbf{e}_{N_1}^{(1) T}$ and $\mathbf{e}_{N_1}^{(1)}$, we find
\begin{multline*}
\mathbf{e}^{(1) T}_{N_1} V \mathbf{e}^{(1)}_{N_1} - \frac{2}{N_1}\mathbf{1}_{N_1}^T V \mathbf{e}^{(1)}_{N_1} + \frac{1}{N_1^2}\mathbf{1}_{N_1}^T V \mathbf{1}_{N_1} \\
= \frac{N_1 - 1}{N_1} \text{ (since $V$ is symmetric).}
\end{multline*}
But since $L_1^T \mathbf{e}^{(1)}_{N_1} = \mathbf{0}$ and $Q_1\mathbf{1}_{N_1} = \mathbf{0}$, we know that both $\mathbf{e}^{(1) T}_{N_1} V \mathbf{e}^{(1)}_{N_1} = 0$ and $\mathbf{1}^{T}_{N_1} V \mathbf{1}_{N_1} = 0$. Thus
\[
\mathbf{1}_{N_1}^T V \mathbf{e}^{(1)}_{N_1} = \frac{1 - N_1}{2} \Rightarrow \mathbf{1}_{N_1}^T L_1 Q_1^T \Sigma_1 Q_1 \mathbf{e}^{(1)}_{N_1} = \frac{1 - N_1}{2},
\]
and so \refe{eqn:lhs2} becomes $LHS = \mathbf{0}$. Thus \refe{eqn:lyapsol} also satisfies \refe{eqn:lyapsolt} and is therefore a solution to the given matrix equation.
\end{IEEEproof}

We can now proceed to state our next main result.

\begin{theorem}\label{theo:notpath}
Suppose $\mathcal{G}_1$ is a connected graph containing only one globally reachable node, and let $\mathcal{G}$ be the graph formed by connecting the globally reachable node in $\mathcal{G}_1$ to an additional node via a directed edge of arbitrary weight. Then the effective resistance between any two nodes in $\mathcal{G}_1$ is equal to the effective resistance between them in $\mathcal{G}$.
\end{theorem}

\begin{IEEEproof}
Let $N_1$, $A_1$, $D_1$ and $L_1$ be the number of nodes, the adjacency matrix, the matrix of node out-degrees and the Laplacian matrix of $\mathcal{G}_1$, respectively. Let $Q_1 \in \mathbb{R}^{(N_1-1)\times N_1}$ satisfy \refe{eqn:propq}. Using $Q_1$, we can compute $\overline{L}_1$ from \refe{eqn:lbar} and since $\mathcal{G}_1$ is connected, we can find matrices $\Sigma_1$ and $X_1$ from \refe{eqn:resdetails}. Without loss of generality, suppose that the globally reachable node in $\mathcal{G}_1$ is node $1$. Since this is the only globally reachable node in $\mathcal{G}_1$, no other nodes can be reached from node $1$ and hence $d_1^\text{\text{\emph{out}}} = 0$. Thus
\begin{equation}\label{eqn:zerodegree}
L_1^T \mathbf{e}^{(1)}_{N_1} = \mathbf{0}.
\end{equation}

Let the additional node in $\mathcal{G}$ be node $N = N_1+1$. We can see that since node $1$ is globally reachable in $\mathcal{G}_1$ and node $N$ is reachable from node $1$, node $N$ is globally reachable in $\mathcal{G}$. Thus $\mathcal{G}$ is connected.

Now, we can write the adjacency matrix of $\mathcal{G}$ as
\[
A = \begin{bmatrix}A_1 & a\mathbf{e}^{(1)}_{N_1}\\\mathbf{0}^T & 0\end{bmatrix},
\]
where $a > 0$ is the weight on edge $(1,N)$ in $\mathcal{G}$. Similarly, we can write the matrix of node out-degrees as
\[
D = \begin{bmatrix}D_1 + a\mathbf{e}^{(1)}_{N_1}\mathbf{e}^{(1) T}_{N_1} & \mathbf{0}\\\mathbf{0}^T & 0\end{bmatrix}.
\]
Utilizing these two expressions, we find that the Laplacian matrix of $\mathcal{G}$ can be written as
\[
L = \begin{bmatrix}L_1 + a\mathbf{e}^{(1)}_{N_1}\mathbf{e}^{(1) T}_{N_1}  & -a\mathbf{e}^{(1)}_{N_1}\\\mathbf{0}^T & 0\end{bmatrix}.
\]

Now, let
\[
Q = \begin{bmatrix}Q_1 & \mathbf{0}\\ \alpha \mathbf{1}_{N_1}^T & -\beta\end{bmatrix},
\]
where $\alpha = \dfrac{1}{\sqrt{N_1\left(N_1+1\right)}}$ and $\beta = \dfrac{N_1}{\sqrt{N_1\left(N_1+1\right)}}$. Then $Q$ satisfies \refe{eqn:propq}. Substituting this matrix into \refe{eqn:lbar} gives us
\[
\overline{L} = \begin{bmatrix}\overline{L}_1 + aQ_1\mathbf{e}^{(1)}_{N_1}\mathbf{e}^{(1) T}_{N_1}Q_1^T & a(\alpha+\beta)Q_1\mathbf{e}_{N_1}^{(1)} \\ \alpha\mathbf{1}_{N_1}^T L_1 Q_1^T + a\alpha \mathbf{e}_{N_1}^{(1) T} Q_1^T & \frac{a}{N_1}\end{bmatrix}.
\]

In order to compute effective resistances in $\mathcal{G}$, we must find the matrix $\Sigma$ which solves \refe{eqn:lyap}. Since we have partitioned $\overline{L}$ into a $2\times 2$ block matrix, we will do the same for $\Sigma$. Let
\[
\Sigma = \begin{bmatrix} S & \mathbf{t}\\ \mathbf{t}^T & u\end{bmatrix},
\]
where $S \in \mathbb{R}^{(N_1 - 1)\times(N_1 - 1)}$ is a symmetric matrix, $\mathbf{t} \in \mathbb{R}^{N_1 - 1}$ and $u \in \mathbb{R}$. Then Lemma \ref{lem:lyapsol} gives a solution to \refe{eqn:lyap} using $\overline{L}$ and our desired form of $\Sigma$. However, since $\mathcal{G}$ is connected, we know that there must be a unique solution. Thus \refe{eqn:lyapsol} from Lemma \ref{lem:lyapsol} is the unique solution to the Lyapunov equation, and in particular,
\[
S = \Sigma_1.
\]

Thus, by Lemma \ref{lem:equalX}, the effective resistance between two nodes in $\mathcal{G}_1$ is equal to the effective resistance between the same two nodes in $\mathcal{G}$.
\end{IEEEproof}

\begin{corollary}\label{cor:removepath}
Suppose $\mathcal{C}_\mathcal{G}(k,j)$ consists of a subgraph $\mathcal{C}^\prime_\mathcal{G}(k,j)$ that is connected via a single edge of arbitrary weight to the leaf node of a directed path. Then the effective resistance between nodes $k$ and $j$ in the graph $\mathcal{G}$ is equal to the effective resistance between nodes $k$ and $j$ in $\mathcal{C}^\prime_\mathcal{G}(k,j)$.
\end{corollary}

\begin{IEEEproof}
This follows by simply applying Theorem \ref{theo:notpath} repeatedly to ``prune'' away the nodes in the directed path.
\end{IEEEproof}

We can see from the graphs shown in Fig.~\ref{fig:pathcounterexample} that this ``pruning'' operation can only be applied to directed edges in general. There may, however, be other graphical structures that also do not affect the effective resistance between two nodes.

\subsection{Effective resistance is a distance-like function and its square root is a metric}\label{subsec:metric}
One useful property of effective resistance in undirected graphs is that it is a metric on the nodes of the graph \cite{Klein93}. This allows effective resistance to substitute for the shortest-path distance in various graphical indices and analyses, as well as offering an alternative interpretation of effective resistance that does not rely on an electrical analogy. Although many of the requirements of a metric follow from the algebraic construction of the effective resistance, the triangle inequality depends on Kirchhoff's laws \cite{Klein93}. Consequently, we shall see that the effective resistance does not satisfy the triangle inequality on general digraphs. Importantly, however, the square root of the effective resistance is a metric. Therefore, if a true metric on digraphs is sought which incorporates information about all connections between two nodes, the square root of the effective resistance is a valid option. In contrast, if the effective resistance had been generalized using the Moore-Penrose generalized inverse instead of our definition, then it would be neither a metric nor would its square root be a metric.

We note that the only difference in the conditions for a metric between a function and its square root lies in the triangle inequality. Furthermore, if a function $d(\cdot,\cdot)$ is a metric, then $\sqrt{d(\cdot,\cdot)}$ is by necessity a metric too.

\begin{theorem}\label{theo:metric}
The square root of the effective resistance is a metric on the nodes of any connected directed graph. That is,
\begin{align}
r_{k,j}	&\geq 0 \;\; \forall \text{ nodes $k$ and $j$,} \label{eqn:rpos}\\
r_{k,j}	&= 0 \; \Leftrightarrow \; k = j,\label{eqn:rdefinite}\\
r_{k,j}	&= r_{j,k} \text{, and} \label{eqn:rsymm}\\
\sqrt{r_{k,\ell}} + \sqrt{r_{\ell,j}}	&\geq \sqrt{r_{k,j}} \;\; \forall \text{ nodes $k$, $j$ and $\ell$.}\label{eqn:rtriangle}
\end{align}
Furthermore, the effective resistance itself is not a metric since it fails to satisfy the triangle inequality.
\end{theorem}

\begin{figure}
\centering
\includegraphics[width=3cm]{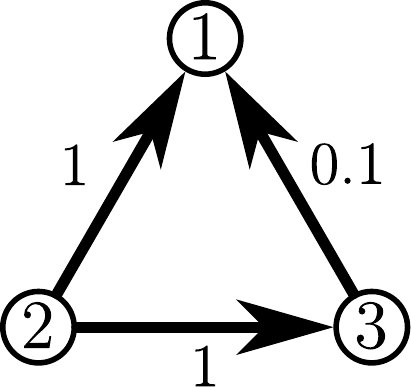}
\caption{A simple 3-node directed graph, $\mathcal{G}^\text{triangle}$,where the triangle inequality fails.}
\label{fig:tricounterexample}
\end{figure}

\begin{IEEEproof}
From \refe{eqn:rdef}, we know that the effective resistance can be computed as 
\[
r_{k,j} = \left(\mathbf{e}^{(k)}_N - \mathbf{e}^{(j)}_N\right)^T X \left(\mathbf{e}^{(k)}_N - \mathbf{e}^{(j)}_N\right),
\]
where $X = 2Q^T \Sigma Q$ and $\Sigma$ is a positive definite matrix. Now, by \refe{eqn:propq}, we know that the matrix $P \mathrel{\mathop :}= \begin{bmatrix}\frac{1}{\sqrt{N}}\mathbf{1}_N & Q^T\end{bmatrix}$ is orthogonal and thus $X$ is similar to $P^T X P = \begin{bmatrix}0 & \mathbf{0}^T\\\mathbf{0} & 2\Sigma\end{bmatrix}$. Hence $X$ has a single $0$ eigenvalue and its remaining eigenvalues are twice those of $\Sigma$. Furthermore, $X\mathbf{1}_N = \mathbf{0}$ since $Q\mathbf{1}_N = \mathbf{0}$. Thus $X$ is positive semi-definite with null space given by the span of $\mathbf{1}_N$.

Since $X$ is positive semi-definite, we can find a matrix $Y \in \mathbb{R}^{N\times N}$ such that $X = Y^T Y$ (e.g. by the Cholesky decomposition or the positive semi-definite square root \cite{Horn85}). This means that we can write effective resistances as $r_{k,j} = \left\|Y \left(\mathbf{e}^{(k)}_N - \mathbf{e}^{(j)}_N\right)\right\|_2^2$, and therefore $\sqrt{r_{k,j}} = \left\|Y \left(\mathbf{e}^{(k)}_N - \mathbf{e}^{(j)}_N\right)\right\|_2$, where $\|\cdot\|_2$ denotes the regular 2-norm on $\mathbb{R}^N$. Therefore, if we associate each node $k$ of $\mathcal{G}$ to the point $\mathbf{p}_k \mathrel{\mathop :}= Y\mathbf{e}^{(k)}_N \in \mathbb{R}^N$, we observe that $\sqrt{r_{k,j}}$ is equal to the Euclidean distance in $\mathbb{R}^N$ between $\mathbf{p}_k$ and $\mathbf{p}_j$. Since $\mathbf{e}^{(k)}_N - \mathbf{e}^{(j)}_N$ is perpendicular to $\mathbf{1}_N$ for any $k \neq j$, $\mathbf{e}^{(k)}_N - \mathbf{e}^{(j)}_N$ is not in the null space of $Y$ and so $\mathbf{p}_k \neq \mathbf{p}_j$ for $k \neq j$. Hence $\sqrt{r_{k,j}}$ is a metric on the nodes of the graph.

Finally, to show that $r_{k,j}$ is not a metric, we consider the 3-node graph shown in Fig.~\ref{fig:tricounterexample}. In this case, we find that $r_{1,3} = 20$, $r_{1,2} = \frac{131}{21} \approx 6.24$ and $r_{2,3} = \frac{37}{7} \approx 5.29$. Thus $r_{1,3} > r_{1,2} + r_{2,3}$ and the triangle inequality does not hold.
\end{IEEEproof}

\begin{remark}
We can observe that if effective resistance was defined using $X = L^\dagger$, then for the graph shown in Fig.~\ref{fig:tricounterexample} the triangle inequality would fail for both $r_{k,j}$ and $\sqrt{r_{k,j}}$. Indeed, the counterexample in the proof of Theorem \ref{theo:metric} also demonstrates why the triangle inequality could be expected to fail for any effective resistance definition that respects edge direction. We can observe that there should be a ``low'' effective resistance between nodes $1$ and $2$ due to the connecting edge with unit weight. Likewise, nodes $2$ and $3$ should have a ``low'' effective resistance between them for the same reason. But node $2$ does not belong to $\mathcal{C}_{\mathcal{G}^\text{triangle}}(1,3)$ and so there should be a ``high'' effective resistance between nodes $1$ and $3$ due to their only connection being an edge with low weight. Thus, the sum of the effective resistances between nodes $1$ and $2$ and between nodes $2$ and $3$ should be lower than the effective resistance between nodes $1$ and $3$.
\end{remark}

\section{Conclusions}\label{sec:conc}
We have generalized the concept of effective resistance to directed graphs in a way that maintains the connection between effective resistances and control-theoretic properties relating to consensus-type dynamics. Despite the algebraic nature of our generalization, 
we have shown that effective resistances in directed graphs bear a fundamental relationship to the structure of the connections between nodes. Moreover, the square root of effective resistance provides a well-defined metric on connected directed graphs, allowing for a notion of distance between nodes, even in cases where neither node is reachable from the other.

Although it may have been tempting to use the Moore-Penrose generalized inverse of a directed graph's Laplacian matrix to define effective resistance, we have shown that not only would this approach ignore the complexity of the derivation of effective resistance for undirected graphs, but also fail to lead to a distance function for directed graphs. Instead, our generalization derives from an analysis of applications of effective resistance in which directed graphs arise naturally. We believe that this approach will allow for the application of this directed version of effective resistance in other situations than those we presented above.

In the companion paper \cite{Young13II}, we demonstrate how to compute effective resistances in certain prototypical classes of graphs and we find cases where effective resistances in directed graphs behave analogously to effective resistances in undirected graphs as well as cases where they behave in unexpected ways.

\FloatBarrier
\bibliographystyle{IEEEtran}
\bibliography{REFabrv,ReferenceList}

\end{document}